\documentclass[a4paper,12pt]{article}

\usepackage{amsfonts}
\usepackage{amssymb}
\usepackage{amsmath}
\usepackage{mathptmx}
\usepackage{color}
\usepackage{cancel}
\usepackage{bbm}
\usepackage{pstricks}
\usepackage{psfrag}
\usepackage{mathrsfs}

\def\ap{'\thinspace}

\def\ns{\hspace{-1mm}}

\newtheorem{lemma}{Lemma}[section]
\newtheorem{theorem}{Theorem}[section]
\newtheorem{remark}{Remark}[section]

\newtheorem{problem}{Problem}[section]
\newtheorem{example}{Example}[section]
\def\be{\begin{equation}}
\def\ee{\end{equation}}
\def\bea{\begin{eqnarray}}
\def\eea{\end{eqnarray}}
\def\beann{\begin{eqnarray*}}
\def\eeann{\end{eqnarray*}}
\def\bsea{\begin{subeqnarray}}
\def\esea{\end{subeqnarray}}
\def\bmat{\left[ \begin{array}}
\def\emat{\end{array} \right]} 
\def\bsmat{\left[ \begin{smallmatrix}}
\def\esmat{\end{smallmatrix} \right]} 

\def\ap{'\thinspace}

\def\proof{\noindent{\bf{\em Proof:}\ \ }}
\def\QED{\mbox{\rule[0pt]{1.5ex}{1.5ex}}}
\def\endproof{\hspace*{\fill}~\QED\par\endtrivlist\unskip}
\def\endex{\hspace*{\fill}~$\square$\par\endtrivlist\unskip}
\newcommand{\real}{{\mathbb{R}}}
\newcommand{\complex}{{\mathbb{C}}}

\def\gR{{\cal R}}

\newcommand{\ima}{\operatorname{im}}
\newcommand{\rank}{\operatorname{rank}}

\newcommand{\diag}{\operatorname{diag}}

\def\tra{{\scalebox{.6}{\thinspace\mbox{T}}}}

\definecolor{Royalblue}{cmyk}{1,0.30,0.2,0.2}

\begin{document}
\begin{titlepage}
\title{\vspace{15mm}
The extended symplectic pencil and the finite-horizon LQ problem with two-sided boundary conditions \thanks{Partially supported by the Italian Ministry for Education and Research (MIUR) under PRIN grant n. 20085FFJ2Z).}\vspace{10mm}}
\author{{\large Augusto Ferrante$^\ddagger$  \quad Lorenzo Ntogramatzidis$^\star$ }\\
       {\small $^\ddagger$Dipartimento di Ingegneria dell\ap Informazione}\\
        {\small    Universit\`a di Padova, via Gradenigo, 6/B -- 35131 Padova, Italy}\\
       {\small     {\tt augusto@dei.unipd.it}} \\ 
       {\small     $^\star$Department of Mathematics and Statistics}\\
       {\small     Curtin University, Perth WA, Australia.}\\
       {\small    {\tt L.Ntogramatzidis@curtin.edu.au}}
 }%
\maketitle
\begin{center}
\begin{minipage}{14.2cm}
\begin{center}
\bf Abstract
\end{center}
This note introduces a new analytic approach to the
solution of a very general class of finite-horizon optimal control problems formulated for discrete-time
systems. This approach provides a parametric expression for 
the optimal control sequences, as well as the corresponding optimal
state trajectories, by exploiting a new decomposition of the so-called extended symplectic pencil. Importantly, the results established in this paper hold under assumptions that are weaker than the ones considered in the literature so far. Indeed, this approach does not require neither the regularity of the symplectic pencil, nor the modulus controllability of the underlying system. In the development of the approach presented in this paper, several ancillary results of independent interest on generalised Riccati equations and on the eigenstructure of the extended symplectic pencil will also be presented.
\end{minipage}
\end{center}
\begin{center}
\begin{minipage}{14.2cm}
\vspace{2mm}
\end{minipage}
\end{center}
\thispagestyle{empty}
\end{titlepage}
\section{Introduction}
\label{secintro}
{T}{his} paper focuses on finite-horizon LQ problems with a generalised cost function and with affine constraints at the end-points. These
problems are not just important {\em per se}. In
fairly recent literature it has been shown that LQ problems are becoming increasingly
useful as building blocks to solve complex
optimisation problems, broken down into two or more LQ subproblems, each one with constraints at
the end-points. A typical example is the $H_2$-optimal tracking problem of previewed reference signals,
which can be solved by splitting the problem into two coupled LQ problems, \cite{Ferrante-N-05}. Finite-horizon LQ problems with constraints at the end-points are also useful in the solution of other optimisation problems, including $H_2$ receding-horizon problems and the minimisation of regulation transient in switching linear plants.

The aim of this paper is to present a method to solve the most general class of finite-horizon linear-quadratic (LQ) optimal control problems in the discrete time with positive semi-definite cost index and affine constraints at the end-points. The approach taken in this paper is based on a procedure for the parameterisation of the set of trajectories generated by the so-called extended symplectic difference equation (ESDE).
The idea of solving finite-horizon LQ problems by exploiting expressions of the trajectories generated by the Hamiltonian system in the continuous time or the ESDE in the discrete time originated in the papers \cite{Ferrante-MN-04}, \cite{Ntogramatzidis-M-1-05} and \cite{Ferrante-N-05} for the continuous time, and in \cite{Ferrante-N-1-04} and \cite{Ferrante-N-06} for the discrete time. In both situations, the expressions parameterising the trajectories of the Hamiltonian system and the symplectic equation hinge on particular solutions of the associated continuous/discrete algebraic Riccati equations and on the solution of the corresponding closed-loop continuous/discrete Lyapunov equation.
While controllability of the given system was required in the first papers \cite{Ferrante-MN-04} \cite{Ferrante-N-1-04}, because both the stabilising and antistabilising solutions of the ARE were involved, in more recent times it has been shown that generalisations of the same technique are possible under the much milder assumption of sign-controllability in the continuous case \cite{Ferrante-N-05} and modulus controllability in the discrete case, see \cite{Ferrante-N-06}. 
In a subsequent paper \cite{EZ} the same problem was considered under the more restrictive assumption of stabilisability. Some problems in the solution presented in \cite{EZ} have been analysed {and corrected} in \cite{Ferrante-N-11}.
The assumptions of sign/modulus-controllability (or stabilisability) were needed in the above-mentioned papers because the solution presented there was based on the existence of a solution of the closed-loop Lyapunov equation.
In the discrete case, the other standing assumption was the regularity of the extended symplectic pencil.
%
%
The goal of this paper is to propose a new approach aimed at overcoming these limitations. More precisely, in this paper a direct method is developed which generalises the technique in \cite{Ferrante-N-1-04} and \cite{Ferrante-N-06} in two directions. First, we do not require the symplectic pencil to be regular, nor to have a spectrum devoid of eigenvalues on the unit circle. As such, with the method proposed in this paper, regular and singular problems can be tackled in a unified manner. Second, unlike the other contributions on this topic, the method presented in this paper does not involve the solvability of the closed-loop Lyapunov equation. Therefore, even the modulus controllability assumption can be dropped. 
The technique presented in this paper only requires a solution of the so-called generalised discrete-time algebraic Riccati equation, which may exist even when the symplectic pencil is not regular (while in this case the standard discrete algebraic Riccati equation cannot be solved). Such solution is used to derive a decomposition of the extended symplectic pencil that yields a natural parameterisation of the solutions of the symplectic difference equation.
A large
number of LQ problems dealt with in the literature by resorting to different -- often iterative -- techniques can be
tackled in a unified framework and in finite, nonrecursive terms, by means of the method developed in this paper. 
{{ For a better description of the features and the generality of our framework, we illustrate all our results in a running example in which the underlying system is not modulus controllable, and the extended symplectic pencil is not regular (so that the methods in previous literature cannot be used).}}


%
%
%
%
\section{Statement of the problem}
\label{stat}
Consider the linear time-invariant discrete-time system governed by the difference equation
\bea
\label{eq1}
x(t+1)=A\,x(t)+B\,u(t), 
\eea
where, for all $t \in \mathbb{N}$, $x(t)\,{\in}\,\mathbb{R}^n$ is the state, $u(t)\,{\in}\,\mathbb{R}^m$
 is the control input, $A\,{\in}\, \mathbb{R}^{n\,{\times}\,n}$ and $B\,{\in}\,\mathbb{R}^{n\,{\times}\,m}$. 
Let $T\,{\in}\,\mathbb{N} \setminus \{0\}$ be the length of the time horizon. 
Let $V_0,V_T\,{\in}\,\real^{q\,{\times}\,n}$ and $v\,{\in}\,\real^q$; consider  
\bea
\label{constraints}
V_0\, x(0)+V_T\,x(T) = v,
\eea
which represents a two-point boundary-value affine constraint on the states at the end-points. With no loss of generality, we can consider $V\triangleq [\, V_0 \;\; V_T \,]$ to be of full row rank.
In the case where $q\,{=}\,0$, the matrices $V_0, V_T, V$ and the vector $v$ are considered to be void: in this case (\ref{constraints}) does not constrain any component of the initial and terminal states.\\[1mm]
 Let $\Pi = \bsmat Q & S \\[1mm] S^{\tra} & R \esmat = \Pi^\tra \ge 0$ be a square $(n+m)$-dimensional matrix 
 with $Q\,{\in}\,\real^{n\,\times\,n}$, $S\,{\in}\,\real^{n\,\times\,m}$ and $R\,{\in}\,\real^{m\,\times\,m}$ (note that we do not assume the non-singularity of $R$). We denote by $\Sigma$ the Popov triple $(A,B,\Pi)$. Finally, let
$H = \bsmat H_1 & H_2 \\[1mm] H_2^{\tra} & H_3 \esmat=H^{\tra} \ge 0$ with $H_1,H_2,H_3\,{\in}\,\real^{n\,\times\,n}$ and $h_0,h_T\,{\in}\,\real^n$. 

\begin{problem}
\label{pro1}
 Find $u(t)$, $t\,{\in}\,\{0,\ldots, T\,{-}\,1\}$ and $x(t)$, $t\,{\in}\,\{0,\ldots, T\}$, minimising 
\bea
 \label{cost}
J(x,u)\ns&\ns \triangleq \ns&\ns 
\sum_{t=0}^{T-1}
\left[ \begin{array}{cc} \!\! x^{\tra}(t)\! &  u^{\tra}(t) \!\!  \end{array} \right]
\Pi \left[ \begin{array}{cc}\! \! x(t) \!\! \\ \! \! u(t) \! \!  \end{array} \right]  +\left[ \begin{array}{cc}  \! \! x^{\tra}(0)-h^{\tra}_0 \!  &  \! x^{\tra}(T)-h_T^{\tra} \!  \!  \end{array} \right]  H  \left[ \begin{array}{c} \! \!   x(0)-h_0 \! \!   \\ \!  \!  x(T)-h_T \!  \!  \end{array} \right], 
\eea
 under the constraints (\ref{eq1}-\ref{constraints}).
 \end{problem}
 The formulation of Problem \ref{pro1} is very general, since the cost index in (\ref{cost}) involves the most general type of positive semidefinite quadratic penalisation on the extreme states, and (\ref{constraints}) represents the most general affine constraint on these states. As particular cases of Problem \ref{pro1} we have:
 \begin{itemize}
 \item the standard case where $x(0)$ is assigned and $x(T)$ is weighted in (\ref{cost}); this case can be recovered from Problem \ref{pro1} by  setting $V_0 = I_n$, $V_T=0$, $h_T\,{=}\,0$, $H_1\,{=}\,H_2\,{=}\,0$; 
  \item the fixed end-point case, where the states at the end-points are sharply assigned; this case can be recovered from Problem \ref{pro1} by setting $V\,{=}\,I_{2n}$ and $H\,{=}\,0$; 
  \item the point-to-point case, where the extreme values of an output $y(t)\,{=}\,C\,x(t)$ are constrained to be equal to two assigned vectors $y_0$ and $y_T$, by taking $V\,{=}\,\textrm{diag}(C,C)$, and $v=\bsmat y_0 \\[1mm] y_T\esmat$. 
  \end{itemize}
  Further non-standard LQ problems that can be useful in practice are particular cases of Problem \ref{pro1}: consider for example an LQ problem in which the states at the end-points $x(0)$ and $x(T)$ are not assigned, but they are constrained to be equal, i.e., $x(0)=x(T)$. This case can be obtained by Problem \ref{pro1} by setting $V_0= I_n$, $V_T = -I_n$ and $v=0$.\\[-0.5cm]
 %

%
\begin{lemma}
\label{lem1}{\bf \cite[Lemma 3]{Ferrante-N-06}}
If $u(t)$ and $x(t)$ are optimal for Problem \ref{pro1}, then $\lambda(t)\,{\in}\,\real^n$, $t\,{\in}\,\{0,\ldots,T\}$ and $\eta\,{\in}\,\real^s$ exist such that $x(t)$, $\lambda(t)$, $u(t)$ and $\eta$ satisfy the set of equations
\bea
&& \hspace{-0.8cm} x(t+1)  =  A\,x(t)+B\,u(t) \qquad t\,{\in}\,\{0,\ldots,T\,{-}\,1\}, \label{eq5} \\
&& \hspace{-0.2cm} V\,\left[ \begin{array}{c} x(0)\\ x(T) \end{array} \right] = v, \label{el2_1}\\
&& \hspace{-0.8cm} \lambda(t) = Q\,x(t)+A^{\top}\,\lambda(t+1)+S\,u(t) \qquad t\,{\in}\,\{0,\ldots,T\,{-}\,1\}, \label{eq5bis} \\
&& \hspace{-0.2cm} \left[ \begin{array}{c} -\lambda(0) \\ \lambda(T) \end{array} \right] = H \left[ \begin{array}{cc} x(0)-h_0 \\ x(T)-h_T \end{array} \right]+V^\top \eta, \label{el3bis}  \\
&& \hspace{-0.8cm} 0 = S^{\top}\,x(t)+B^{\top}\,\lambda(t+1)+R\,u(t) \qquad t\,{\in}\,\{0,\ldots,T\,{-}\,1\}. \label{eq5bisbis} 
%
\eea
 Conversely, if equations (\ref{eq5}-\ref{eq5bisbis}) admit solutions $x(t)$, $u(t)$, $\lambda(t)$, $\eta$, then $x(t)$, $u(t)$ minimise $J(x,u)$ subject to the constraints (\ref{eq1}-\ref{constraints}).
\end{lemma}
The variables $\lambda(t)$ in (\ref{eq5}-\ref{eq5bisbis}) represent 
the Lagrange multipliers associated with the constraint (\ref{eq1}), \cite{Lancaster-95,Ionescu-OW-99}, while the variable $\eta\,{\in}\,\mathbb{R}^s$ is the Lagrange multiplier vector associated with (\ref{constraints}). 


\section{The generalised Riccati equation and the extended symplectic system}
\label{ESDE}
Since in the present setting we are not assuming that $R$ is positive definite, (\ref{eq5bisbis}) cannot be solved in $u(t)$ to obtain a set of $2n$ equations in $x(t)$ and $\lambda(t)$.
A convenient form in which (\ref{eq5}), (\ref{eq5bis}) and (\ref{eq5bisbis}) can be written, that does not require inversion of $R$, is the  descriptor form
\bea
\label{eq6}
M\,p(t+1)=N\,p(t) \qquad t\,{\in}\,\{0,\ldots,T\,{-}\,1\},
\eea
where \\[-1cm]
\[
M\triangleq\left[ \begin{array}{ccc}
 I_n & O & O \\
O & -A^{\tra} & O \\
O & -B^{\tra} & O
\end{array} \right], \quad N\triangleq\left[ \begin{array}{ccc}
A & O & B \\
Q & -I_n & S \\
S^{\tra} & O & R
\end{array} \right], \quad 
p(t)\triangleq\left[ \begin{array}{c}
x(t)\\
\lambda(t) \\
u(t)
\end{array} \right].
\]
The matrix pencil $N-z\,M$ is known as the
{\em extended symplectic pencil}, \cite{Lancaster-95,Ionescu-OW-99}, herein denoted concisely by ESP($\Sigma$). In this paper we do not make the assumption of regularity of this pencil.

%
We now show how a solution of a generalised discrete algebraic Riccati equation can be used to obtain a decomposition of ESP($\Sigma$) that can be used to solve Problem \ref{pro1}. In particular, we will exploit the solutions of the following constrained matrix equation
\bea
X \ns&\ns = \ns&\ns  A^\tra\,X\,A-(A^\tra\,X\,B+S)(R+B^\tra\,X\,B)^{\dagger}(B^\tra\,X\,A+S^\tra)+Q, \label{gdare} \\
\ns&\ns \ns&\ns \qquad \ker (R+B^\tra\,X\,B) \subseteq \ker (A^\tra\,X\,B+S), \label{kercond}
\eea
where the matrix inverse that appears in the standard discrete algebraic Riccati equation (DARE) has been replaced by the Moore-Penrose pseudo-inverse. Eq. (\ref{gdare}) is known in the literature as the {\em generalised discrete-time algebraic Riccati equation} GDARE($\Sigma$), \cite{Stoorvogel-S-98,Ionescu-O-96}. GDARE($\Sigma$) with the additional constraint given by (\ref{kercond}) is sometimes referred to as {\em constrained generalised discrete-time algebraic Riccati equation} CGDARE($\Sigma$). 
Clearly (\ref{gdare}) constitutes a generalisation of the classic DARE($\Sigma$), in the sense that any solution of DARE($\Sigma$) is also a solution of GDARE($\Sigma$) -- and therefore also of CGDARE($\Sigma$) -- but the {\em vice-versa} is not true in general.

 We now introduce some notation that will be used throughout the paper. First, to any matrix $X=X^\tra \in \real^{n \times n}$  we associate the following matrices:
\bea
S_X  &  \triangleq  &  A^\tra X\, B \! + \! S,\quad
R_X   \triangleq   R \! + \! B^\tra X B, \quad G_X   \triangleq  I_m-R_X^\dagger R_X, 
\label{defgx}
 \\
K_X  &  \triangleq  &  (R+B^\tra\,X\,B)^\dagger (B^\tra\,X\,A+S^\tra)=R_X^\dagger S_X^\tra, \quad A_X   \triangleq  A-B\,K_X.  \label{KX}
\eea 
The term $R_X^\dagger R_X$ is the orthogonal projector that projects onto $\ima R_X^\dagger=\ima R_X$ so that $G_X$ is  the orthogonal projector that projects onto $\ker R_X$. Hence,
 $\ker R_X=\ima G_X$.

Since as aforementioned the Popov matrix $\Pi$ is assumed to be symmetric and positive semidefinite, we can consider a standard factorisation of the form $\Pi=\bsmat C^\tra \\[1mm] D^\tra\esmat \bsmat  C \,&\, D \esmat$, where $Q=C^\tra C$, $S=C^\tra D$ and $R=D^\tra D$. 

 \begin{example}
 \label{ex0}
 Consider the following Popov triple, which will be used as a running example throughout the paper:
 \beann
A=\left[ \begin{array}{cc} 1 & 1 \\ 0 & 1 \end{array} \right], \quad B=\left[ \begin{array}{ccc} 2 & 0 \\ 1 & 1 \end{array} \right], \quad Q=\left[ \begin{array}{ccc}  0 & 0 \\ 0 & 1  \end{array} \right], \quad S=\left[ \begin{array}{ccc} 0 & 0 \\ 0 & 0 \end{array} \right], \quad R=\left[ \begin{array}{ccc}  0 & 0 \\ 0 & 0  \end{array} \right].
\eeann
The Popov matrix $\Pi$ in this case can be factored with $C=[\,0\;\;1\,]$ and $D=[\,0\;\;0\,]$. The extended symplectic pencil in this case is not regular. 
As such, DARE($\Sigma$) in this case does not admit solutions. On the other hand, in this case CGDARE($\Sigma$) admits the solution $X=\diag\{0,1\}$, that can be computed by resorting to the algorithm proposed in \cite{Ferrante-04}. In this case, $R_X=R+B^\tra\,X\,B=\left[ \begin{smallmatrix} 1 & 1 \\[1mm] 1 & 1 \end{smallmatrix} \right]$, and the corresponding closed-loop matrix is $A_X=\diag\{1,0\}$. Observe that the spectrum of $A_X$ is not unmixed. 
 \endex
\end{example}

 Let $C_X \triangleq C-D\,R_X^\dagger S_X^\tra$, and let $\gR_X$ denote the reachable subspace associated with the pair $(A,B\,G_X)$, in symbols $\gR_X \triangleq \ima [\, B G_X \;\;\; A_X B G_X  \;\;\; A_X^2 B G_X  \;\;\; \ldots  \;\;\; A_X^{n-1} B G_X \,]$.

The following results were proved in \cite[Lemma 4.1, Lemma 4.2, Theorem 4.2]{Ferrante-N-12-sub}. 

\begin{lemma}
\label{lemmaa}
Let $X=X^\tra$ be a solution of CGDARE($\Sigma$). Then,
\begin{description}
\item{\bf (i)} $\gR_X \subseteq\ker C_X$;
\item{\bf (ii)} $\ker R_X =  \ker (XB) \cap \ker R$;
\item{\bf (iii)} $\gR_X$ is the reachability subspace on the output-nulling subspace $\ker X$.
\end{description}
\end{lemma}

We also have the following results, see \cite[Theorems 4.3-4.4]{Ferrante-N-12-sub}.

\begin{lemma}
\label{lemmab}
Let $X$ and $Y$ be two solutions of CGDARE($\Sigma$). Let $A_X$ and $A_Y$ be the corresponding closed-loop matrices. Then,
\begin{description}
\item{\bf (i)} $\ker R_{X}=\ker R_{Y}$;
\item{\bf (ii)} $\gR_{X}=\gR_{Y}$;
\item{\bf (iii)} {$ A_X|_{\gR_{X}}= A_Y|_{\gR_{Y}}$. } 
\end{description}
\end{lemma}

The following result adapts \cite[Lemma 2.5]{Ferrante-W-07} to the case when the matrix pencil $N-z\,M$ may be singular.

\begin{lemma}
\label{the0}
{\mbox Let $X=X^\tra$ be a solution of CGDARE($\Sigma$). Then, $U_X,V_X \in \real^{2n\!+\!m}$ exist such that}
\be\label{tiblokdec}
U_X\,(N-z\,M)\,V_X=\left[ \begin{array}{ccc}
A_X-z\,I_n & O & B \\
O & I_n-z\,A_X^\tra & O \\
O & -z\,B^\tra &  R_X \end{array} \right].
\ee
\end{lemma}
\proof By direct computation we find 
\beann
U_X (N\!-\!z M) V_X\!=\!\!\left[ \begin{array}{ccc}
 \!\!\! A_X-z\,I_n \!\! & \!\! O \!\! & \!\! B \!\!\! \\
 \!\!\! \Xi_{21} \!\! & \!\! I_n-z\,A_X^\tra \!\! & \!\! \Xi_{23} \!\!\! \\
 \!\!\! \Xi_{23}^\tra \!\! & \!\! -z\,B^\tra \!\! & \!\!  R_X  \!\!\! \end{array} \right]\,\;\textrm{with}\;\,\;
U_X \triangleq \left[ \begin{array}{ccc} 
 \!\!\! I_n  \!\! & \!\! O \!\! & \!\! O \!\!\!\! \\
 \!\!\! A^\tra_X X \!\! & \!\! I_n \!\! & \!\! -K_X^\tra \!\!\!\! \\
 \!\!\! B^\tra X \!\! & \!\! O \!\! & \!\! I_m  \!\!\!\! \end{array} \right]\!\!, \;
V_X \triangleq\left[ \begin{array}{ccc} 
\! \!\!\! I_n \!\! & \!\!\! O \!\! & \!\! O \!\!\! \\ 
 \!\!\!\! X \!\! & \!\!\! -I_n \!\! & \!\! O \!\!\! \\ 
 \!\!\!\! -K_X \!\! & \!\!\! O \!\! & \!\! I_m \!\!\! \end{array} \right]\!\!.
\eeann
The term $\Xi_{21}$ is given by
\beann
\Xi_{21} \ns&\ns = \ns&\ns A_X^\tra X\,A-A_X^\tra X\,B\,K_X+Q-X-S\,K_X-K_X^\tra S^\tra +K_X^\tra R \, K_X -z\,(A^\tra\,X-A_X^\tra X+K_X^\tra B^\tra X).
\eeann
The term multiplying $z$ is zero since $A_X=A-B\,K_X$. Moreover, since GDARE($\Sigma$) can be written as $X=A^\tra X\,A-S_X\,K_X+Q$ we find
$\Xi_{21} = K_X^\tra (R_X\,K_X-S_X^\tra)=S_X\,R_X^\dagger R_X\,R_X^\dagger S_X^\tra-S_X\,R_X^\dagger S_X^\tra=0$.
Finally, $\Xi_{23}=A^\tra X B-z\,X\,B - K_X^\tra B^\tra X B+S+z\,X\,B-K_X^\tra R= S_X\,G_X$. 
In view of (\ref{kercond}), we have $S_X\,G_X=0$, so that (\ref{tiblokdec}) holds. \endproof

If $X$ is a solution of CGDARE($\Sigma$), from the triangular structure in (\ref{tiblokdec}) we have
\bea
\label{det}
\det(N-z\,M)=\det(A_X-z\,I_n)\cdot \det(I_n-z\,A_X^\tra) \cdot \det R_X.
\eea

When $R_X$ is non-singular {(i.e. $X$ is a solution of DARE($\Sigma$))}, the dynamics represented by this matrix pencil are decomposed into a part governed by the generalised eigenstructure of $A_X-z\,I_n$, a part governed by the finite generalised eigenstructure of $I_n-z\,A_X^\tra$, and a part which corresponds to the dynamics of the eigenvalues at infinity. When $X$ is a solution of DARE($\Sigma$), the generalised eigenvalues\footnote{Recall that a generalised eigenvalue of a matrix pencil $N-z\,M$ is a value of $z \in \complex$ for which the rank of the matrix pencil $N-z\,M$ is lower than its normal rank.} of $N-M\,z$ are given by the eigenvalues of $A_X$, the reciprocal of the non-zero eigenvalues of $A_X$, and a generalised eigenvalues at infinity whose algebraic  multiplicity is equal to $m$ plus the algebraic multiplicity of the eigenvalue of $A_X$ at the origin, and we have
\bea
\label{division}
\sigma(N-z\,M)=\sigma(A_X-z\,I_n) \cup \sigma \left( \left[ \begin{array}{cc} I_n-z\,A_X^\tra & O \\
 -z\,B^\tra &  R_X\end{array} \right] \right). 
 \eea

When the matrix $R_X$ is singular, (\ref{det}) still holds but provides no information as in this case $\det R_X=0$, while (\ref{division}) is no longer true. We show this fact with a simple example.

 \begin{example}
 \label{ex1}
 Consider Example \ref{ex0}. 
Matrix $X=\diag\{0,1\}$ is a solution of CGDARE($\Sigma$), and the corresponding closed-loop matrix is $A_X=\diag\{1,0\}$. From Lemma \ref{the0} we find
\beann
U_X\,(N-z\,M)\,V_X=\left[ \begin{array}{cc|cc|cc}
 1 - z &  0 &     0 &  0 & 2 & 0 \\
     0 & -z &     0 &  0 & 1 & 1 \\
     \hline
     0 &  0 & 1-z &  0 & 0 & 0\\
     0 &  0 &     0 & 1 & 0 & 0\\
     \hline
     0 &  0 &   -2\,z &  -z & 1 & 1\\
     0 &  0 &     0 &  -z & 1 & 1 \end{array} \right],
     \eeann
{whose normal rank (which coincides with that of  $N-z\,M$) is easily seen to be} equal to $5$. The eigenvalues of $A_X$ are $0$ and $1$. However, it is not true that $z=1$ is a generalised eigenvalue of $N-z\,M$. In fact, a direct check shows that the rank of $N-M$ is equal to $5$.\footnote{{We warn that the routine {\tt eig.m} of the software MATLAB$^{\textrm{\tiny{\textregistered}}}$ (version 7.11.0.584(R2010b)) in this case fails to provide the right answer. It indeed returns $1$ as a generalised eigenvalue of the pencil $N-z\,M$.} 
  }
  \endex
 \end{example}
 \ \\[-1cm]
 
Consider a change of coordinates in the input space $\real^m$ induced by {the $m \times m$ orthogonal matrix $T=[\, T_{1} \;\; T_{2} \,]$ where $\ima T_{1}=\ima R_X$ and $\ima T_{2}=\ima G_{X}=\ker R_X$.} From Lemma \ref{lemmab}, $T$ is independent of the solution $X$ of CGDARE($\Sigma$). Thus 
$T^\tra R_X\,T= \diag\{R_{X,0},O\}$, where $R_{X,0}$ is invertible. Its dimension is denoted by $m_1$. Consider the block matrix $\hat{T}\triangleq\diag(I_n, I_n, T)$. Defining  the matrices $B_1 \triangleq B\,T_{1}$ and $B_2 \triangleq B\,T_{2}$ 
 we get
 \beann
 \hat{T}^\tra \left( U_X\,(N-z\,M)\,V_X \right) \hat{T}= 
 \left[ \begin{array}{cccc}
 A_X-z\,I_n & O & B_1 & B_2 \\
O & I_n-z\,A_X^\tra & O & O\\
O & -z\,B_1^\tra &  R_{X,0} & O \\
O & -z\,B_2^\tra &  O & O 
 \end{array} \right].
 \eeann
In view of $\ker R_X=\ima G_X$, we get $\ima B_2= \ima (B\,G_X)$. Matrix $B_1$ has $m_1$ columns. Let $m_2 \triangleq m - m_1$ be the number of columns of $B_2$. Let us take
 $U=[\,U_1\;\;U_2\,]$ such that $\ima U_1$ is the reachable subspace associated with the pair $(A_X,B_2)$, which coincides with the subspace $\gR_X$. We have
 \bea
 \label{kalman}
 U^{-1}\,A_X\,U=\bmat{cc}  A_{X,11} & A_{X,12}  \\
 O & A_{X,22}\emat, \qquad U^{-1}\,B_2=\bmat{cc}  B_{21} \\ O  \emat, \qquad U^{-1}\,B_1=\bmat{cc}  B_{11} \\ B_{12} \emat.
 \eea
  Let $\hat{U}=\diag\{U,U,I_{m_1},I_{m_2}\}$. 
 Let $r$ denote the size of $\gR_X$. Defining the two unimodular matrices
 \beann
 \Omega_1 \triangleq \left[ \begin{array}{cccccc}
 \!  I_r  \! & \! O  \! & \! O  \! & \! O  \! & \! O  \! & \! O  \! \\
 \!   O  \! & \! O  \! & \! I_r  \! & \! O  \! & \! O  \! & \! O  \! \\
  \!   O  \! & \! O  \! & \! O  \! & \! O  \! & \! O  \! & \! I_{m_2}  \! \\
  \!    O  \! & \! I_{n-r}  \! & \! O  \! & \! O  \! & \! O  \! & \! O  \! \\
  \!     O  \! & \! O  \! & \! O  \! & \! I_{n-r}  \! & \! O  \! & \! O  \! \\
  \!      O  \! & \! O  \! & \! O  \! & \! O  \! & \! I_{m_1}  \! & \! O 
       \end{array} \right] \quad \textrm{and} \quad
       \Omega_2 \triangleq \left[ \begin{array}{cccccc}
 \!  I_r  \! & \! O  \! & \! O  \! & \! O  \! & \! O  \! & \! O  \! \\
  \!  O  \! & \! O  \! & \! O  \! & \! I_{n-r}  \! & \! O  \! & \! O  \! \\
  \!   O  \! & \! O  \! & \! I_r  \! & \! O  \! & \! O  \! & \! O  \! \\
   \!   O  \! & \! O  \! & \! O  \! & \! O  \! & \! I_{n-r}  \! & \! O  \! \\
   \!    O  \! & \! O  \! & \! O  \! & \! O  \! & \! O  \! & \! I_{m_1}  \! \\
      \!  O  \! & \! I_{m_2}  \! & \! O  \! & \! O  \! & \! O  \! & \! O 
       \end{array} \right],\quad \textrm{we get}
       \eeann
  \bea
&&\hspace{-1.7cm} P(z) \triangleq \Omega_1\,\hat{U}^{-1}\,\hat{T}^\tra \left( U_X\,(N-z\,M)\,V_X \right) \hat{T}\,\hat{U}\, \Omega_2= \nonumber \\
&& =
       \left[ \begin{array}{cc|c|ccc}
 A_{X,11}-z\,I_{r} & B_{21} & O & A_{X,12} & O & B_{11} \\
 \hline
 O & O & I_{r}-z\,A_{X,11}^\tra & O & O & O \\
 O & O & -z\,B_{21}^\tra & O & O & O \\
 \hline
 O & O & O & A_{X,22}-z\,I_{n-r}& O & B_{12} \\
 O & O & -z\,A_{X,12}^\tra & O & I_{n-r}\,-z\,A_{X,22}^\tra & O \\
 O & O & -z\,B_{11}^\tra & O & -z\,B_{12}^\tra & R_{X,0} 
  \end{array} \right]. \label{fc}
  \eea
Since the pair $(A_{X,11}, B_{21})$ is reachable by construction, all the $r$ rows of the submatrix $[\,A_{X,11}-z\,I_{r} \;\; B_{21} \,]$ are linearly independent for every $z \in \complex \cup \{\infty\}$. This also means that of the $r+m_2$ columns of $[\,A_{X,11}-z\,I_{r} \;\; B_{21} \,]$, only $r$ are linearly independent, and this gives rise to the presence of a null-space of $P(z)$ whose dimension $m_2$ is independent of $z \in \complex \cup \{\infty\}$. We obtain\footnote{Let $\Xi=\bsmat \Xi_{11} & \Xi_{12} \\[1mm] O & \Xi_{22} \esmat$. Observe that if either $\Xi_{11}$ is full row-rank or  $\Xi_{22}$ is full column-rank, then $\rank\, \Xi=\rank\, \Xi_{11}+\rank \,\Xi_{22}$.}
\beann
\rank P(z)=r+ \rank\left[ \begin{array}{c|ccc}
    I_{r}-z\,A_{X,11}^\tra & O & O & O \\
 -z\,B_{21}^\tra & O & O & O \\
 \hline
  O & A_{X,22}-z\,I_{n-r}& O & B_{12} \\
  -z\,A_{X,12}^\tra & O & I_{n-r}-z\,A_{X,22}^\tra & O \\
  -z\,B_{11}^\tra & O & -z\,B_{12}^\tra & R_{X,0} 
  \end{array} \right].
  \eeann
  Now, consider the rank of $\left[ \begin{smallmatrix} 
  I_{r}-z\,A_{X,11}^\tra  \\ -z\,B_{21}^\tra \end{smallmatrix} \right]$. Again, since the pair $(A_{X,11}, B_{21})$ is reachable, this rank is constant and equal to $r$ for every $z \in \complex \cup \{\infty\}$. Thus, 
  \beann
  \rank P(z)=2\,r+ \rank P_1(z), \quad \textrm{where} \quad P_1(z) \triangleq \left[ \begin{array}{ccc}
A_{X,22}-z\,I_{n-r}& O & B_{12} \\
 O & I_{n-r}-z\,A_{X,22}^\tra & O \\
 O & -z\,B_{12}^\tra & R_{X,0} 
  \end{array} \right].
  \eeann
 Since $\det P_1(z)=\det(A_{X,22}-z\,I_{n-r})\cdot \det(I_{n-r}-z\,A_{X,22}^\tra) \cdot \det R_{X,0}$, a value $z \in \complex$ can be found for which $\det P_1(z) \neq 0$. Hence, the normal rank of $P_1(z)$ is equal to $2\,(n-r)+m_1$, and therefore the normal rank of $P(z)$ is
 $2\,r+2\,(n-r)+m_1=2\,n+m_1$. The generalised eigenvalues of the pencil $P(z)$ are the values $z \in \complex \cup \{\infty\}$ for which the rank of $P_1(z)$ is smaller than its normal rank $2\,(n-r)+m_1$. These values are the eigenvalues of $A_{X,22}$ plus their reciprocals, included possibly the eigenvalue at infinity, whose multiplicity {| be it algebraic or geometric | is, in general, {\em not} given by the sum of $m_1$ plus the multiplicity of  the eigenvalue in zero of $A_{X,22}$.
In fact, the multiplicity of the eigenvalue at infinity is the 
multiplicity of the zero eigenvalue of}
\beann
 P_{\infty}\triangleq \left[ \begin{array}{ccc} I_{n-r} & O & O \\ O & A_{X,22}^\tra & O \\ O & B_{12}^\tra & O  \end{array} \right].
\eeann
  {If $A_{X,22}$ is non-singular, the last $m_1$ columns give rise to an eigenvalue at infinity whose multiplicity (algebraic and geometric) is exactly equal to $m_1$, since in this case the dimension of the null-space of $P_{\infty}$ is equal to $m_1$.   
However, if $A_{X,22}$ is singular, the algebraic (geometric) multiplicity of the zero eigenvalue of  $P_{\infty}$ is equal the sum of $m_1$ plus the algebraic (geometric) multiplicity of the  eigenvalue in zero of $A_{X,22}^\tra$
that is non-observable for the pair $(A_{X,22}^\tra,B_{12}^\tra)$
and these are indeed the multiplicities of the eigenvalue at infinity of the pencil. }\\[0cm]

  From these considerations, it turns out that, unlike the regular case, not all the eigenvalues of $A_X$ appear as generalised eigenvalues of ESP($\Sigma$). In particular, the eigenvalues of $A_X$ restricted to $\gR_{X}$ do not appear as generalised eigenvalues, whereas the eigenvalues of the map induced by $A_X$ in the quotient space $\real^n /  \gR_{X}$ -- along with the reciprocals of those that are different from zero -- are generalised eigenvalues of ESP($\Sigma$).

 \begin{example}
 Consider Example \ref{ex0}. Using the solution $X=\diag \{0,1\}$ of CDARE($\Sigma$) we get $\ker R_X=\bsmat -1 \\[1mm] 1 \esmat$ and $\ima R_X=\bsmat 1 \\[1mm] 1 \esmat$.  By taking $T=\left[ \begin{smallmatrix} 1 & -1 \\[1mm] 1 & 1 \end{smallmatrix} \right]$ we obtained
$T^\tra R_X\,T=
\diag\{4,0\}$.
 Hence, in this case $m_1=m_2=1$. We partition $B\,T$ as $B\,T=\left[ \begin{smallmatrix} 2 & -2 \\[1mm] 2 & 0 \end{smallmatrix} \right]$, so that $B_1=\left[ \begin{smallmatrix} 2  \\[1mm] 2  \end{smallmatrix} \right]$ and $B_2=\left[ \begin{smallmatrix} -2 \\[1mm]  0 \end{smallmatrix} \right]$.
The normal rank of ESP($\Sigma$) is equal to $2\,n+m_1=5$.
The generalised eigenvalues of $N-z\,M$ are given by the uncontrollable eigenvalues of the pair $(A_X,B_2)=\left(\left[ \begin{smallmatrix} 1 & 0 \\[1mm] 0 & 0 \end{smallmatrix} \right],\left[ \begin{smallmatrix} -2 \\[1mm] 0 \end{smallmatrix} \right]\right)$ plus their reciprocals. Therefore, ESP($\Sigma$) has a generalised eigenvalue at the origin. Since $A_{X,22}=0$ and $B_{12}=2$, it also has an eigenvalue at infinity with multiplicities equal to the multiplicities of the zero eigenvalue of $\left[ \begin{smallmatrix} 0 & 0 \\[1mm] 2 & 0 \end{smallmatrix} \right]$. 
By writing this pencil in the form given by (\ref{fc}), we get 
 \beann
 \hat{T}^\tra \left( U_X\,(N-z\,M)\,V_X \right) \hat{T}= 
 \left[ \begin{array}{cc|c|ccc}
  1-z & -2 & 0 & 0 & 0 & 2 \\
  \hline
    0 & 0 & 1-z & 0 & 0 & 0 \\
      0 & 0 & 2\,z & 0 & 0 & 0 \\
        \hline
        0 & 0 & 0 & -z & 0 & 2 \\
          0 & 0 & 0 & 0 & 1 & 0 \\
            0 & 0 & -2\,z & 0 & -2\,z & 4 \end{array} \right],
 \eeann
 from which we see that zero is indeed the only finite generalised eigenvalue of ESP($\Sigma$).
 \endex
 \end{example}

\section{Solution of the LQ problem}
In the basis constructed in the previous section, (\ref{eq6}) can be written for $t \in \{0,\ldots,T-1\}$ as
\bea
x_1(t+1) & = & A_{X,11}\, x_1(t)+B_{21}\, u_1(t) + A_{X,12}\, x_{2}(t)+B_{11}\, u_2(t), \label{eqa} \\
\lambda_1(t) & = & A_{X,11}^\tra \,\lambda_1(t+1), \label{eqb} \\
0 & = & -B_{21}^\tra\, \lambda_1(t+1), \label{eqc} \\
x_2(t+1) & = & A_{X,22}\, x_2(t)+B_{12}\, u_2(t), \label{eqd} \\
\lambda_2(t) & = & A_{X,22}^\tra\, \lambda_2(t+1)+A_{X,12}^\tra\, \lambda_1(t+1), \label{eqe} \\
u_2(t) &=& R_{X,0}^{-1} B_{12}^\tra\, \lambda_2(t+1)+R_{X,0}^{-1}\, B_{11}^\tra\, \lambda_1(t+1). \label{eqf} 
\eea
Since by construction the pair $(A_{X,11},B_{21})$ is reachable, $\ker \bsmat A_{X,11}^\tra \\[1mm] B_{21}^\tra \esmat=\{0\}$, which means (\ref{eqb}-\ref{eqc}) yield $\lambda_1(t)=0$ for all $t \in \{0,\ldots,T-1\}$. This implies that (\ref{eqe}-\ref{eqf}) can be simplified as
\bea
\lambda_2(t) & = & A_{X,22}^\tra\, \lambda_2(t+1), \label{eqe1} \\
u_2(t) &=& R_{X,0}^{-1} \,B_{12}^\tra \,\lambda_2(t+1). \label{eqf1} 
\eea
It is clear at this point that we can parameterise all the trajectories generated by the difference equations (\ref{eqd}), (\ref{eqe1}) and (\ref{eqf1}) in terms of $x_2(0)$ and $\lambda_2(T)$. Indeed, (\ref{eqe1}) leads to
\bea
\label{eqe2}
\lambda_2(t)=(A_{X,22}^\tra)^{T-t} \, \lambda_2(T) \qquad \forall\,t \in \{0,\ldots,T\}.
\eea
This expression can be plugged into (\ref{eqf1}) and leads
\bea
\label{eqf2}
u_2(t)=R_{X,0}^{-1}\, B_{12}^\tra\, (A_{X,22}^\tra)^{T-t-1} \, \lambda_2(T).
\eea
Plugging (\ref{eqe2}) and (\ref{eqf2}) into (\ref{eqd}) gives
\bea
\label{eqd2}
x_2(t)=A_{X,22}^t x_2(0)+\sum_{j=0}^{t-1} A_{X,22}^{t-j-1} B_{12} R_{X,0}^{-1} B_{12}^\tra (A_{X,22}^\tra)^{T-j-1} \, \lambda_2(T).
\eea
It is worth observing that 
\bea
\label{sum}
x_2(T)=A_{X,22}^T x_2(0)-P \, \lambda_2(T), \quad \textrm{where} \quad P \triangleq \sum_{j=0}^{T-1} A_{X,22}^{T-j-1} B_{12} R_{X,0}^{-1} B_{12}^\tra (A_{X,22}^\tra)^{T-j-1}.
\eea
It is easy to see that matrix $P$ can be re-written as $P=\sum_{j=0}^{T-1} A_{X,22}^{j} B_{12} R_{X,0}^{-1} B_{12}^\tra (A_{X,22}^\tra)^{j}$.
Therefore, $P$ satisfies the discrete Lyapunov equation
\beann
P=A_{X,22}\,P\,A_{X,22}^\tra-A_{X,22}^T\,B_{12} R_{X,0}^{-1} B_{12}^\tra\,(A_{X,22}^\top)^T+B_{12} R_{X,0}^{-1} B_{12}^\tra.
\eeann
If $A_{X,22}$ has unmixed spectrum, this equation can be used to determine $P$ instead of computing the sum in (\ref{sum}). At this point we can solve (\ref{eqa}), which can be written as
\bea
x_1(t+1) & = & A_{X,11}\, x_1(t)+B_{21}\, u_1(t) + \xi(t),
\eea
where $\xi(t)=A_{X,12}\,x_2(t)+B_{11}\,u_2(t)$. Using (\ref{eqd2}) and (\ref{eqf2}) we find
\beann
\xi(t)\!=\!A_{X,12} A_{X,22}^t x_2(0)\!+\!\!\!\left(\!\!B_{11} R_{X,0}^{-1} B_{12}^\tra (A_{X,22}^\tra)^{T\!-\!t\!-\!1}\!+\!A_{X,12}\!\sum_{j=0}^{t-1} \!A_{X,22}^{t\!-\!j\!-\!1} B_{12} R_{X,0}^{-1} B_{12}^\tra (A_{X,22}^\tra)^{T-j-1}\! \!\right)\! \lambda_2(T).
\eeann
{
Let $R_1=[B_{21}\mid A_{X,11}\,B_{21}\mid A_{X,11}^2\,B_{21}\mid\dots\mid A_{X,11}^{T-1}B_{21}]$ and
$R_2=[I\mid A_{X,11}\mid A_{X,11}^2\mid \dots\mid A_{X,11}^{T-1}]$. Then, we can write $x_1(T)=A_{X,11}^T x_1(0)+R_2\,\Xi+R_1\,U_1$ where $\Xi \triangleq \bsmat \xi(T-1) \\[-3mm] \vdots \\ \xi(0) \esmat$ and $U_1 \triangleq \bsmat u_1(T-1) \\[-3mm] \vdots \\ u_1(0) \esmat$. We assume that $T$ is greater than the controllability index of the pair $(A_{X,11},B_{21})$. All the solutions of this equation are parameterised by
\bea
\label{eqa3}
U_1=R_1^\dagger\left(x_1(T)-A_{X,11}^T x_1(0)-R_2\,\Xi\right)+(I-R_1^\dagger\,R_1)\,v_1.
\eea
where $v_1$ is arbitrary.
}

\subsection{Boundary conditions}
Consider the change of coordinates given by the matrix $U=[U_1\;\;U_2]$, where $\ima U_1$ is the reachable subspace of the pair $(A_X,B\,G_X)$. Let $\bsmat x_1(t) \\[1mm] x_2(t) \esmat=U^{-1}\,x(t)$ be the coordinates of the state in the basis induced by $U$, partitioned conformably with $U$. The state, co-state and transversality equations can be written again as in (\ref{eq5}), (\ref{eq5bis}) and (\ref{eq5bisbis}), where $A$, $B$, $Q$, $S$, $V$, $H$, $h_0$ and $h_T$ are replaced by $U^{-1}\,A\,U$, $U^{-1}\,B$, $U^\tra\,Q\,U$, $U^\tra\,S$, $V\,\bsmat U & O \\[1mm] O & U \esmat$, $\bsmat U & O \\[1mm] O & U \esmat^\tra H\,\bsmat U & O \\[1mm] O & U \esmat$, $U^{-1}h_0$ and $U^{-1}h_T$, respectively. 
We can now write (\ref{el2_1}) and (\ref{el3bis}) with respect to this basis. 
We can eliminate the multiplier $\eta$ from (\ref{el3bis}) by premultiplying both sides of this equation by a basis $K_V$ of $\ker V$:
\bea
\label{el3bisbis}
K_V^\tra H \bmat{c} x(0) \\ x(T) \emat + K_V^\tra \bmat{cc} I & O \\ O & -I \emat  \bmat{c} \lambda(0) \\ \lambda(T) \emat = K_V^\tra H \bmat{c} h_0 \\ h_T \emat.
\eea
In this way, (\ref{el2_1}) and (\ref{el3bisbis}) can be written together as a set of $2\,n$ linear equations in $x(0)$, $x(T)$, $\lambda(0)$ and $\lambda(T)$. However, in (\ref{sum}) the component $x_2(T)$ is expressed as a linear function of $x_2(0)$ and $\lambda_2(T)$, and $\lambda_2(0)$ can be expressed as a linear function in $\lambda_2(T)$ by (\ref{eqe2}). Finally we know that $\lambda_1(t)$  must be identically zero, so that $\lambda_1(0)=\lambda_1(T)=0$. {Therefore, in this basis (\ref{el2_1}) and (\ref{el3bis}) can be expressed as a single linear equation of the form
\bea
\label{lineq}
F\,x=g, \quad\textrm{where}\quad  x=\bmat{cccc} x_1^\tra(0) & x_1^\tra(T) & x_2^\tra(0) & \lambda_2^\tra(T) \emat^\tra.
\eea
We have just proved the following result.
\begin{theorem}
Problem \ref{pro1} admits solutions if and only if (\ref{lineq}) does. For any solution $x=[\, x_1^\tra(0) \;\;\; x_1^\tra(T)  \;\;\;  x_2^\tra(0)  \;\;\;  \lambda_2^\tra(T) \,]^\tra$ we get an optimal initial state $x(0)=\bsmat x_1(0) \\[1mm] x_2(0) \esmat$ and a class of optimal controls parameterised by (\ref{eqf2}) and (\ref{eqa3}). The solutions obtained in this way are all the solutions of Problem \ref{pro1}.
\end{theorem}
}

 \begin{example}
 \label{ex2}
 Consider a finite-horizon LQ problem in the time interval $\{0,\ldots,T\}$, involving the matrices given in Example \ref{ex0}.
The initial and final states are constrained to be equal, i.e., $x(0)=x(T)$. Let $H=I_{2\,n}$, $h_0=\bsmat h_1 \\[1mm] h_2 \esmat$ and $h_T=0$.
As aforementioned, $X=\diag\{0,1\}$ is a solution of CGDARE($\Sigma$), 
 leading to $A_X=\diag\{1,0\}$. 
By taking $T=\left[ \begin{smallmatrix} 1 & -1 \\[1mm] 1 & 1 \end{smallmatrix} \right]$, we obtained
 $T^\tra R_X\,T=\diag\{4,0\}$, so that $R_{0,X}=4$.  Recall that $B\,T=\left[ \begin{smallmatrix} 2 & -2 \\[1mm] 2 & 0 \end{smallmatrix} \right]$, so that $B_1=\left[ \begin{smallmatrix} 2  \\[1mm] 2  \end{smallmatrix} \right]$ and $B_2=\left[ \begin{smallmatrix} -2 \\[1mm]  0 \end{smallmatrix} \right]$.
 Therefore, the reachable subspace of the pair $(A_X,B_2)$ is $\ima \bsmat 1 \\[1mm] 0 \esmat$, which means this system is already in the desired basis. Thus, $A_{X,11}=1$, $A_{X,12}=A_{X,22}=0$, $B_{11}=B_{12}=2$ and $B_{21}=-2$.
 In this case,  (\ref{eqd}), (\ref{eqe1}) and (\ref{eqf1}) become
\beann
x_2(t+1) = B_{12}\, u_2(t), \qquad \lambda_2(t) = 0 \cdot \lambda_2(t+1), \qquad
u_2(t)= R_{X,0}^{-1} B_{12}^\tra \lambda_2(t+1). \\[-1.4cm]
\eeann

This implies that
\beann
\lambda_2(t)=\left\{ \begin{array}{ll} 0 & t \in \{0,\ldots,T-1\} \\
\lambda_2(T) & t = T,
\end{array} \right. \quad \rightarrow \quad u_2(t)=\left\{ \begin{array}{ll} 0 & t \in \{0,\ldots,T-2\} \\
R_{X,0}^{-1} B_{12}^\tra \,\lambda_2(T) & t = T-1,
\end{array} \right.
\eeann
 which give\\[-1cm]
 \beann
x_2(t)=\left\{ \begin{array}{ll} x_2(0) & t=0 \\
0 & t \in \{1,\ldots,T-1\} \\
B_{12}^\tra R_{X,0}^{-1} B_{12}^\tra \lambda_2(T) = \lambda_2(T) & t = T.
\end{array} \right.
\eeann
 In this basis, (\ref{el2_1}) gives rise to $x_{1}(0) = x_1(T)$ and $x_2(0) = x_2(T)=\lambda_2(T)$, 
which are linear in $x_1(T)$ and $\lambda_2(T)$,
while (\ref{el3bisbis}) can be written as
$x_1(0)+x_1(T)=h_1$ and $x_2(0)+x_2(T)+\lambda_2(0)-\lambda_2(T)=h_2$. Since $\lambda_2(0)=0$ and $x_2(T)=\lambda_2(T)$, the latter can be written as $x_2(0)=h_2$. 
Therefore, the boundary conditions can be written in the form (\ref{lineq}):
\beann
\bmat{cccc} 1 & 0 & -1 & 0 \\ 0 & 1 & 0 & -1\\ 1 & 0 & 1 & 0 \\ 0 & 1 & 0 & 0 \emat \bmat{c} x_1(0) \\ x_2(0) \\ x_1(T) \\ \lambda_2(T) \emat=
\bmat{c} 0 \\ 0 \\ h_1 \\ h_2 \emat.
\eeann
This linear equation admits only the solution $x_1(0)=x_1(T)=h_1/2$ and $x_2(0)=\lambda_2(T)=h_2$.
Now we can compute the optimal control law. First, $u_2(t)$ is zero for all $t \in \{0,\ldots,T-2\}$ and $u_2(T-1)=R_{X,0}^{-1} B_{12}^\tra \,\lambda_2(T)=h_2/2$. In order to compute $u_1$, we write (\ref{eqa}) as
\bea
\label{zeta}
x_1(t+1) = 1\cdot x_1(t)-2\, u_1(t) + \xi(t).
\eea
The term $\xi(t)$ in this case is equal to zero for all $t \in \{0,\ldots,T-2\}$ and $\xi(T-1)=B_{11}\,R_{0,X}^{-1}\,B_{12}^\tra \,\lambda_2(T)=\lambda_2(T)= h_2$.
 We can write (\ref{eqa3}) explicitly as\\[-1cm]
\beann
x_1(T)=x_1(0)+\bmat{cccccc} I & A_{X,11} & A_{X,11}^2 & \ldots & A_{X,11}^{T-1} \emat\!\! \bmat{c} h_2 \\ 0 \\ \vdots \\[-3mm] 0 \emat+
\underbrace{\bmat{cccccc} -2 & -2 & \ldots & -2 \emat}_{T} \!\! \bmat{c} u_1(T-1) \\ u_1(T-2) \\ \vdots \\[-3mm] u_1(0) \emat
\eeann
which gives \\[-1cm]
\beann
 \bmat{c} u_1(T-1) \\ u_1(T-2) \\ \vdots \\[-3mm] u_1(0) \emat= \frac{h_2}{2\,T}\bmat{c}1  \\ 1 \\ \vdots \\[-3mm] 1 \emat+
 \bmat{cccccc} 
 1-T & 0 &  \ldots & 0  \\
  1 & 2-T &  \ldots & 0 \\
  \vdots & \vdots& \ddots & \vdots \\[-3mm]
  1 & 1 &  \ldots & -1 \\
   1 & 1 &  \ldots & 1\emat\,v,
   \eeann
   where $v$ is arbitrary and represents the degree of freedom in the control $u_1$. 
 \endex
\end{example}

\begin{remark}
So far, we have not considered the problem of existence of solutions for Problem \ref{pro1}. In general, the existence of a state trajectory $x(t)$ satisfying the constraints (\ref{eq1}-\ref{constraints}) for some $u(t)$ is not ensured, since we have not assumed reachability on (\ref{eq1}). A necessary and sufficient condition for the existence of optimal solutions is that there exist state and input trajectories satisfying (\ref{eq1}-\ref{constraints}) {(feasible solutions)}.  
In fact, since the optimal control problem formulated in Section \ref{stat} involves a finite number of variables -- precisely, $L=m\cdot T$ for the control plus  $n$ for the initial state -- Problem \ref{pro1} can be restated as a quadratic static optimization problem in these $L+n$ variables with linear constraints. Thus, a solution to Problem \ref{pro1} exists if and only if a feasible solution -- i.e., a state and input functions satisfying both (\ref{eq1}) and (\ref{constraints}) -- exists.
\end{remark}
\begin{remark}
The approach presented in this paper can successfully tackle even more general LQ problems, where the performance index is not necessarily positive semidefinite. E.g., consider 
\[
J(x,u)=
\sum_{t=0}^{T-1}
\left[ \begin{array}{cc} \!\! x^{\tra}(t)\! &  u^{\tra}(t) \!\!  \end{array} \right]
\Pi \left[ \begin{array}{cc}\! \! x(t) \!\! \\ \! \! u(t) \! \!  \end{array} \right]  + x^{\tra}(T)\,H\, x(T)+2\,\zeta^\tra\,x(T).
\]
Although all the variational analysis remains unaffected, the presence of the term $2\zeta \, x(T)$ deserves some considerations. Indeed, this linear term may cause the divergence to $-\infty$ of the cost index in correspondence to a sequence of  admissible controls so that, even in the presence of feasible solutions, the optimal control may fail to exist.\footnote{Consider for example the case where $A$, $B$ and $Q$ are the $2 \times 2$ identity matrices, while $R$, $S$ and $H$ are the zero matrices and $\zeta=[\, 1 \;\; 1 \,]^\top$. For this system, the LQ problem in one step (i.e., $T=1$) has no solution; in fact, the control
 $u(0)=-x(0) - m\,\zeta$ yields a value of the cost which goes to $-\infty$ as the parameter $m$ goes to $+\infty$.}
In this case, the linear equation representing the boundary conditions is infeasible.
Two simple {\em a priori} sufficient conditions for the existence of the optimal control and hence for the solvability of the two-point boundary-value problem are the following:
\begin{enumerate}
\item $\ker H\subseteq \ker \zeta^\top$. Under this condition, the cost on the final state (and hence  the overall cost index)  is bounded from below.
Indeed, such a cost may be rewritten as a constant plus a positive semi-definite quadratic form $(x(T)-\bar{x})^\top H (x(T)-\bar{x})$ in the difference between $x(T)$ and a suitable ``target state" $\bar{x}$.
 In this case the solution of the problem indeed exists. 
\item $R>0$. In this case the current cost increases quadratically with the norm of the control input with the largest norm and, in the best situation, decreases linearly with the same norm.
Thus the search for the optimal control input can be restricted to a compact set in 
$\real^{m\times T}$ and hence the optimal solution does exist.
\end{enumerate}
\end{remark}

\section{Conclusions}
In this note, we studied the discrete-time finite-horizon LQ problem with the most general type of positive semidefinite cost function and with affine constraints at the end-points. We derived an analytic approach, based on a special decomposition of the extended symplectic pencil, that generalises several contributions that have appeared in the literature in the last few years on this problem. Indeed, this approach does not require regularity of the extended symplectic pencil, nor the modulus controllability of the underlying system. Due to its generality, the proposed technique can be used to efficiently tackle complex optimisation problems of wide interest, including the $H_2$-optimal rejection/tracking of previewed signals, receding-horizon optimal control problems, and the minimisation of regulation transients for plants subject to large parameter jumps.

\end{document}